\documentclass[10pt]{article}
\usepackage{graphicx}
\usepackage{psfrag}
\usepackage{amsmath} 
\usepackage{amsfonts}
\usepackage{amssymb}
\newtheorem{theorem}{Theorem}

\newtheorem{lemma}[theorem]{Lemma}

\newenvironment{proof}[1][Proof]{\textbf{#1.} }{\ \rule{0.5em}{0.5em}}

\newcommand{\DD}{{\mathbb D}}
\newcommand{\oDD}{{\overline{\mathbb D}}}
\newcommand{\CC}{{\mathbb C}}
\newcommand{\RR}{{\mathbb R}}
\newcommand{\Dir}{{\mathcal D}}
\renewcommand{\a}{\alpha}
\renewcommand{\o}{\omega}
\renewcommand{\b}{\beta}

\newcommand{\pp}{{p^{\prime}}}
\newcommand{\vp}{{\varphi}}
\newcommand{\oT}{{\overline{T}}}
\renewcommand{\SS}{{\mathcal S}}
\newcommand{\oS}{{\overline{\mathcal S}}}
\newcommand{\DoS}{{\overline{ S}}}

\begin{document}

\title{The characterization of the Carleson
measures for analytic Besov spaces: a simple proof.}
\author{N. Arcozzi\thanks{Partially supported by the COFIN project Analisi Armonica,
funded by the Italian Minister for Research}\\Dipartimento di Matematica\\Universit\`a di Bologna\\40127 Bologna, ITALY
\and R. Rochberg\thanks{This material is based upon work supported by the National
Science Foundation under Grant No. 0400962.}\\Department of Mathematics\\Washington University\\St. Louis, MO 63130, U.S.A.
\and E. Sawyer\thanks{This material based upon work supported by the National
Science and Engineering Council of Canada.}\\Department of Mathematics \& Statistics\\McMaster University\\Hamilton, Ontario, L8S 4K1, CANADA}
\maketitle

\section{Introduction and summary}
In this note we return on the characterization of the Carleson measures for the class of weighted analytic Besov spaces $B_p(\rho)$ considered in \cite{ARS1}, Theorem 1. Since our first paper on this topic, we have found shorter or better proofs for some of the intermediate results, which were used in proving various extensions and generalizations of the characterization theorem, \cite{A}\cite{ARS2}\cite{ARS3}. Here we assemble these new arguments to give the characterization theorem a self-contained, short and simple proof. Recall that a positive measure $\mu$ on $\DD$ is Carleson for $B_p(\rho)$ if
$$
\int_{\DD}|f|^pd\mu\le C(\mu)^p\|f\|_{B_p(\rho)}^p.
$$
Here and below, $C(\mu)$ is a generic constant depending on $\mu$, $\rho$ and $p$ alone. The basic idea in proving Theorem 1 in \cite{ARS1} was that, for $f$ in $B_p(\rho)$, the quantity
$$
(1-|z|^2)|f^\prime(z)|=\Delta\psi
$$
can be assumed to be essentially constant, at least on a metaphorical level, on Whitney boxes in the unit disc $\DD$. Here $\Delta$ is the backward difference operator on the dyadic tree $T$ of the Whitney boxes, $\Delta\psi(\a)=\psi(\a)-\psi(\a^{-1})$, and its inverse operator is the discrete Hardy operator $I$: $I\varphi(\a)=\sum_{o\le\b\le\a}\varphi(\b)$, where $o$ is the \it root \rm of $T$. The metaphor can be made precise in such a way that the problem of characterizing the Carleson measures is completely reduced to the characterization of the Carleson measures for $B^T_p(\rho)$, a discrete version of $B_p(\rho)$ \cite{A}\cite{ARS2}. Such characterization was proved in \cite{ARS1} by means of a direct, but technical good-$\lambda$ inequality. 
Here we reproduce an alternative, short proof from \cite{ARS3}, consisting of a simple interpolation argument.

In this note we make two observations which lead to an improvement of the Carleson measure theorem. The first is that the discrete function $\psi=I\Delta\psi$ is directly related to the radial variation of the Besov function $f$. The second is that the characterization of the Carleson measures for $B^T_p(\rho)$, with the new proof, allows us to consider measures $\mu$ living on $\oDD$, the closure of $\DD$. This leads to the proof that a measure $\mu$ on $\oDD$ is Carleson for $B_p(\rho)$ if and only if the Carleson inequality holds with $f$ replaced by its radial variation $V(f)$,
$$
\int_{\oDD}V(f)^pd\mu\le C(\mu)^p\|f\|_{B_p(\rho)}^p.
$$
The result for the radial variation is related with Beurling's Theorem on exceptional sets for the Dirichlet spaces \cite{B}. We will return on this subject in future work.

There is a vast literature on Carleson measures for weighted analytic Besov spaces. See, for instance, \cite{Car}\cite{Ste}\cite{KS}\cite{V}\cite{CV}\cite{CO}. We single out Carleson's pioneering paper \cite{Car}, dealing with the Hardy space ($H^2=B_2(1-|z|^2)$ can be seen as a weighted Dirichlet space, which is not included in the family of the spaces $B_p(\rho)$ considered in this note), from which the Carleson measures took their name; and Stegenga's article \cite{Ste}, which was the first dealing with \it bona fide \rm Dirichlet spaces. Stegenga realized that the Carleson measures for the Dirichlet space $\Dir=B_2(1)$ could not, unlike the Hardy case, be characterized by a simple "one-box condition" of Carleson type. He gave a capacitary condition to be checked on finite unions of Carleson boxes. Later on, non-capacitary, one-box characterizations were given and this note presents one of them.

Here is an outline of the paper. in \S 2 we state the main theorem and we give some preliminary results, in \S 3
we develop the needed discrete material, in \S 4 we finish proof of the main theorem by a duality argument. The logical structure of the proof and most of the technical details are the same in the general case and in the Dirichlet case: the reader might set $p=2$ and $\rho\equiv1$ without missing the essential features of the proof.

\smallskip

{\small The first author would like to thank the organizers of the Thessaloniki conference for their kind hospitality.}

\section{The main theorem}

Let $1<p<\infty$ and let $\rho$ be a pointwise positive, measurable weight on $\DD$. The {\it analytic Besov space $B_p(\rho)$} contains the functions $f$ which are holomorphic in $\DD$ such that
\begin{equation}
\label{normabp}
\|f\|_{B_p(\rho)}=
\left\{\left|(1-|z|^2)f^\prime(z)\right|^p\rho(z)\frac{dA(z)}{(1-|z|^2)^2}\right\}^{1/p}+|f(0)|<
+\infty.
\end{equation}
Here, $dA(z)=\frac{dxdy}{\pi}$ is normalized area measure. See \cite{Zhu} for the general theory of these spaces in the unweighted case.
Note that $B_2(1)=\Dir$ is the Dirichlet space, for which we have an inner product that we denote by $<\cdot,\cdot>_{\Dir}$.
Throughout this note we assume that the weight $\rho$ is 
 {\it $p$-admissible}, i.e. that 
\begin{enumerate}
\item[(dual)] the dual of $B_p(\rho)$ under the duality pairing $<\cdot,\cdot>_{\Dir}$ is $B_{p^\prime}(\rho^{1-p^\prime})$;
\item[(reg)] there are constants $0<c<1$ and $C>0$ such that $\rho(z)\le\rho(w)$ whenever 
$$
\left|\frac{z-w}{1-\overline{z}w}\right|\le c.
$$
\end{enumerate}
Condition (reg) says that $\rho$ changes by a bounded factor in balls of fixed hyperbolic radius. Recall the hyperbolic metric is given by
$$
ds^2=\frac{|dz|^2}{(1-|z|^2)^2}.
$$ 
Although we do not make explicit use of it, the hyperbolic geometry is the right setting for thinking of analytic Besov spaces.

A characterization of the $p$-admissible weights in terms of a boundary $A_p$ condition follows easily from Bekoll\'e's Theorem \cite{B} on Bergman projections; see \cite{ARS1} \S4 and \S5, for a discussion. We mention that $\rho(z)=(1-|z|^2)^a$, $a\ge0$, is $p$-admissible for $0\le a<1$.
 
A positive Borel measure $\mu$ on $\oDD$ is a {\it Carleson measure} for $B_p(\rho)$ if for all $f$ in $B_p(\rho)$ the radial limits
$$
f(e^{it})=\lim_{r\to1}f(re^{it})
$$
exists $\mu$-a.e for $e^{it}$ in $\partial\DD$ and if
\begin{equation}
\label{carleson}
\left\{\int_\oDD|f^\prime(z)|^p d\mu(z)\right\}^{1/p}\le C(\mu) \|f\|_{B_p(\rho)}.
\end{equation}
In literature, generally $\mu$ is defined as a measure on $\DD$. We can however extend any such $\mu$ by $\tilde{\mu}(E)=\mu(E\cap\DD)$, so that the requirement on the radial limits is vacuously fulfilled and (\ref{carleson}) holds with $\tilde{\mu}$ instead of $\mu$. 

The {\it radial variation} of a function $f$, holomorphic in $\DD$, is the function $V(f):\oDD\to[0,+\infty]$,
$$
V(f)(Re^{it})=\int_0^R|f^\prime(re^{it})|dr.
$$
We can now state the characterization theorem for the Carleson measures. 
\begin{theorem}
\label{main}
Let $1<p<\infty$ and let $\rho$ be a $p$-admissible weight. Let $\mu$ be a positive Borel measure on $\oDD$. Then, the following are equivalent.
\begin{enumerate}
\item[(Car)] $\mu$ is a Carleson measure for $B_p(\rho)$.
\item[(VarCar)] $\mu$ satisfies the Carleson inequality for the radial variation:
\begin{equation}
\label{variation}
\left\{\int_{\oDD}|V(f)(z)|^pd\mu(z)\right\}^{1/p}\le C(\mu)\|f\|_{B_p(\rho)}.
\end{equation}
\item[(TreeCar)] $\mu$ is a Carleson measure for the discrete Besov space $B_p^T(\rho)$.
\item[(TC)] $\mu$ satisfies the \it tree condition \rm
\begin{equation}
\label{TC}
\sum_{\b\ge\a}\mu(\DoS(\b))^\pp(\b)\rho^{1-\pp}(\b)\le C(\mu)\mu(\DoS(\a)).
\end{equation}
\end{enumerate}
\end{theorem}
Throughout this note, $\pp$ denotes the exponent conjugate to $p$, $p^{-1}+\pp^{-1}=1$.
See below in this section for an explanation of the terminology used in (TreeCar) and (TC).

The proof of Theorem \ref{main} is divided into several steps. At the end of this section we show that (TreeCar)$\implies$(VarCar)$\implies$(Car). In \S 4 we will show that (Car)$\implies$(TreeCar). In particular, the problem of characterizing the Carleson measures for $B_p(\rho)$ is reduced to a discrete one, the equivalence of (TreeCar) and (TC), which will be proved in \S3. 

\smallskip

Consider, now, a dyadic Whitney decomposition of
$\mathbb{D}$. Namely, for integer $n\geq0,\ 1\leq m\leq2^{n}$, let
\[
D_{n,m}=\left\{  z\in\mathbb{D}\colon2^{-n-1}\leq1-|z|\leq2^{-n}%
,\ \frac{m-1}{2^n} \le {\frac{\arg(z)}{2\pi}} \le \frac{m}{2^n}  \right\}  .
\]
Each Whitney box is approximatively a hyperbolic ball of unitary radius. Namely,
$$
\left|\frac{z-w}{1-\overline{z}w}\right|\le c<1
$$ 
for $z,w$ in $D_{n,m}$, with $c$ independent of $n$ and $m$.

It is natural to consider the Whitney squares as indexed by the vertices of a
dyadic tree, $T$. Thus the vertices of $T$ are
\begin{equation}
\left\{  \alpha|\alpha=(n,m),\ n\geq0\ \mbox{and}\ 1\leq m\leq2^{n}%
,\ m,n\in\mathbb{N}\right\} \label{eqdeltaalpha}%
\end{equation}
and we say that there is an edge between $(n,m),\ (n^{\prime},m^{\prime})$ if
$D_{(n,m)}$ and $D_{(n^{\prime},m^{\prime})}$ share an arc of a
circle. The \it root \rm of $T$ is, by definition, $o=(0,1)$. Here and throughout we
will abuse notation and, when convenient, identify the vertices of a such a
tree with the sets for which they are indices: $\Delta_\a=\a$. Thus, there are two edges having $(0,1)$ as
endpoint, each other box being the endpoint of exactly three edges.
Given points $\a,\b$ in $T$, the \it geodesics \rm joining $\a$ and $\b$, $[\a,\b]$, is the set of vertices one has to cross to go from $\a$ to $\b$, following the edges of $T$. The distance between $\a$ and $\b$ in $T$ is $d(\a,\b)=\sharp([\a,\b])-1$.

The choice of the root induces a partial ordering on $T$: $\a\le\b$ if $\a\in[o,\b]$. For $\a$ in $T$, we define
$$
S(\a)=\cup_{\b\ge\a}\b,\ \DoS(\a)=\overline{S(\a)},\ \mbox{is\ the\ closure\ of\ }\a\ \mbox{in\ }\DD,\ \partial S(\a)=\DoS(\a)-S(\a).
$$
The immediate predecessor of $\a\ne o$ is the unique point $\a^{-1}=\b$ on $T$ such that $\b\in[o,\a]$ and $d(\a,\b)=1$.
Given a function $\psi: T\to\CC$, we consider the \it backward difference operator \rm
$$
\Delta\psi(\a)=\begin{cases}
\psi(o) & \mbox{if}\ \a=o\crcr
\psi(\a)-\psi(\a^{-1}) & \mbox{if}\ \a\ne o
               \end{cases}
$$
and its inverse, the \it discrete Hardy operator \rm
$$
I\psi(\a)=\sum_{o\le \b\le\a}\psi(\b).
$$
We have $\Delta\circ I=I\circ\Delta=Id,$ the identity.
Since an admissible weight satisfies (reg), we can identify it with a weight $\tilde{\rho}$ on $T$. For instance, by letting $\tilde{\rho}(\a)=\frac{1}{|\a|}\int_\a\rho(z)dA(z)$. With slight abuse of language, we call $\tilde{\rho}=\rho$, using the same name for the weight in $\DD$ and in $T$.   
 The \it discrete Besov  space \rm $B_p^T(\rho)$ is contains those $\psi:T\to\CC$ such that
 $$
 \|\psi\|_{B_p^T(\rho)}^p=\sum_{\a\in T}|\Delta\psi(\a)|^p\rho(\a)<\infty.
 $$
Given a positive measure $\mu$ on $\DD$, we associate to it a measure on $T$ by the obvious formula $\mu(\a)=\int_\a d\mu(z)$. The measure $\mu$ is \it a Carleson measure for \rm $B_p^T(\rho)$ if the inequality
$$
\sum_{\a\in T}|\psi(\a)|^p\mu(\a)\le C(\mu)\|\psi\|_{B_p^T(\rho)}^p
$$
holds.

However, we want to consider measures $\mu$ with support in $\oDD$, hence we need consider boundary measures on the tree. The reader who is just interested in Carleson measures supported on $\DD$ can skip the following paragraph. 

The \it boundary \rm $\partial T$ of $T$ consists of the infinite tree geodesics $\omega$ leaving $o$, $\omega=[\omega_0=o,\omega_1,\dots,\omega_n,\cdots)$, with $\omega_j$ in $T$. Here, for each $j\ge1$, $\omega_{j-1}$ and $\omega_j$ are endpoints of an edge. We extend the partial ordering to 
the closure of $T$, $\overline{T}=T\cup\partial T$, setting $\omega_j\le\omega$ when $\omega$ lies in $\partial T$ and $\omega_j$ is an element of $\omega$.  We give $\overline{T}$ the topology having as basis the sets
$$
\oS(\a)=\{\omega\in\oT:\ \omega\ge\a\},\ \a\in T.
$$
Also, consider the sets
$$
\SS(\a)=\{\omega\in T:\ \omega\ge\a\},\ \partial \SS(\a)=\{\omega\in \partial T:\ \omega\ge\a\},\ \a\in T.
$$
For each $\a\in T$, viewed as a Whitney box, let $I(\a)=\overline{S(\a)}\cap\partial\DD$ be the corresponding boundary arc. 
The map $\Phi:\omega=[\omega_0,\omega_1,\cdots)\mapsto\cap_jI(\omega_j)$ is a correspondence from $\partial T$ onto $\partial\DD$, which is $1-1$ but for a countable set and which sends Borel sets into Borel sets. 
Given a measure $\nu$ on $\partial\DD$, let $\tilde{\nu}=\Phi^*\nu$ be the pull-back of $\nu$ through $\Phi$. Since a measure $\mu$ on $\oDD$ can be decomposed as $\mu=\mu_\DD+\mu_{\partial\DD}$, with $\mu_{\partial\DD}$ supported on $\partial\DD$ and $\mu_\DD(\partial\DD)=0$, we can pull-back any positive Borel measure on $\oDD$ to a Borel measure $\tilde{\mu}$ on $\oT$. 
With abuse of language, we set $\tilde{\mu}=\mu$.
We extend the definition above by saying that the measure $\mu$ supported on $\oT$ is {\it a Carleson measure for} $B_p^T(\rho)$ if 
$$
\int_{\oT}|\psi(\omega)|^pd\mu(\omega)\le C(\mu)\|\psi\|_{B_p^T(\rho)}^p.
$$

The inequality can be rephrased in terms of the operator $I$,
\begin{equation}
\label{paris}
\int_{\oT}|I\varphi(\omega)|^pd\mu(\omega)\le C(\mu)\|\varphi\|_{L^p(T,\rho)}^p,
\end{equation}
where $I$ extends to boundary points $\omega$ in $\partial T$ by
$$
I\varphi(\omega)=\sum_{\a\in\omega}\varphi(\a)
$$
and where $\|\varphi\|_{L^p(T,\rho)}^p=\left(\sum_{\a\in T}|\varphi(\a)|^p\rho(\a)\right)^{1/p}$.

\smallskip

We continue with some elementary estimates.

To each funtion $f$ in $B_p(\rho)$ we associate a function $\vp:T\to[0,+\infty)$ as follows. If $\a\in T$ and $\a\ne o$,
$$
\Delta\vp(\a)=(1-|z(\a)|)|f^\prime(z(\a))|,
$$
where $z(\a)\in\overline{\a}$ and
$$
|f^\prime(z(\a))|=\sup_{w\in\a}|f^\prime(w)|.
$$
We also set $\vp(o)=|f(0)|$.
\begin{lemma}
\label{norm}
$\|f\|_{B_p(\rho)}\sim\|\vp\|_{B_p^T(\rho)}.$
\end{lemma}
\begin{proof}[Proof of the Lemma]
\begin{eqnarray*}
\|f\|_{B_p(\rho)}^p&\sim&|f(0)|^p+\sum_{\a\in T}\int_\a |(1-|z|^2)f^\prime(z)|^p\rho(z)\frac{dA(z)}{(1-|z|^2)^2}\crcr
&\le& |f(0)|^p+C\sum_{\a\in T}\int_\a \frac{dA(z)}{(1-|z|^2)^2}\cdot|(1-|z(\a)|)f^\prime(z(\a))|^p\rho(z(\a))\crcr
&\sim&\sum_{\a\in T}(\Delta\a)^p\rho(\a)\crcr
&=&\|\vp\|_{B_p^T(\rho)}^p.
\end{eqnarray*}
In the other direction, consider for each $\a\in T$ the disc $B(\a)$, a hyperbolic disc having hyperbolic radius independent of $\a$, such that $\a$ is contained in $B(\a)$ and such that the hyperbolic distance between $\partial\a$ and $\partial B(\a)$ is bounded below by some positive $c$ (or, which is the same, such that the Euclidean distance between the two boundaries is bounded below by $c(1-|\a|)$).  By the Mean Value Property and Jensen's inequality we have, if $\a\ne o$,
\begin{eqnarray*}
(\Delta\vp(\a))^p&=&((1-|z(\a)|)|f^\prime(z(\a))|)^p\crcr
&=&(1-|z(\a)|)^p\left|\frac{\pi}{|B(\a)|}\int_{B(\a)}f^{\prime}(z)dA(z)\right|^p\crcr
&\le&\pi^p(1-|z(\a)|)^p\frac{1}{|B(\a)|}\int_{B(\a)}\left|f^{\prime}(z)\right|^p dA(z)\crcr
&\sim&\int_{B(\a)}\left|(1-|z|^2)f^{\prime}(z)\right|^p \frac{dA(z)}{(1-|z|^2)^2}.
\end{eqnarray*} 
Summing over $\a$ and taking into account the fact that each point in $\DD$ belongs to a universally bounded number of discs $B(\a)$, we have
\begin{eqnarray*}
\|\vp\|_{B_p^T(\rho)}^p&\le&|f(0)|^p+C\sum_{\a\in T}\rho(\a)
\int_{B(\a)}\left|(1-|z|^2)f^{\prime}(z)\right|^p \frac{dA(z)}{(1-|z|^2)^2}\crcr
&\sim&\|f\|_{B_p(\rho)}^p.
\end{eqnarray*}
\end{proof}

The function $\varphi$ can be extended to a function $\varphi:\oT\to[0,+\infty]$, defining
$$
\varphi(\omega)=\sum_{j=1}^\infty\Delta\varphi(\omega_j)
$$
when $\omega$ belongs to $\partial T$.
\begin{lemma}
\label{tollo}
If $\o\in T$ and $Re^{it}\in\o$ or if $\o\in\partial T$ and $Re^{it}=e^{it}\equiv\o$, then
$$
V(f)(Re^{it})\le C\vp(\o).
$$
\end{lemma}
\begin{proof}
Let $\omega\in \overline{T}$ and let $[\o,o]$ be the geodesic between $o$ and $\o$ in $\overline{T}$. We have:
\begin{eqnarray*}
V(f)(Re^{it})&\le&\sum_{}|f(0)|+\int_{z(\o_{j-1})}^{z(\o_j)}|f^\prime(z)|d|z|\crcr
&\le &C\left(|f(0)|+\sum_{}(1-|z(\o_j)|)|f^\prime(z(\o_j))|\right)\crcr
&=&\sum_{}\Delta\vp(\o_j)=C\vp(\o).
\end{eqnarray*}
\end{proof}

From Lemmas \ref{norm} and \ref{tollo} we have that (TreeCar)$\implies$(VarCar), 
$$
\int_{\oDD}V(f)^pd\mu\le \int_{\oT}\varphi^pd\mu\le C(\mu)\|\varphi\|_{B_p^T(\rho)}^p\le C(\mu)\|f\|_{B_p(\rho)}.
$$
On the other hand, (VarCar) trivially implies (Car).
\section{Carleson measures on trees.}
We prove that (TreeCar) is equivalent to (TC). 
We rewrite (TreeCar) in its dual form,
\begin{equation}
\label{milano}
\sum_{\a\in T}\left(\int_{\oS(\a)}gd\mu\right)^\pp\rho^{1-\pp}(\a)\le C(\mu)\int_{\oT}g^\pp d\mu\ \mbox{for}\ g\ge0.
\end{equation}
In fact, (TreeCar) is equivalent to the boundedness of $I:L^p(T,\rho)\to L^p(\oT,\mu)$ (see (\ref{paris})), and this is in turn equivalent to the boundedness of the formal adjoint $I^*_\mu$ of $I$,
$$
I^*_\mu:L^\pp(\oT,\mu)\to L^\pp(T,\rho^{1-\pp}),
$$
where duals are taken with respect to the $L^2(T,1)$ and the $L^2(\oT,\mu)$ inner products, respectively,
$$
<f,I^*_\mu g>_{L^2(T,1)}=<If,g>_{L^2(\oT,\mu)}.
$$
Below, we use the obvious fact that $\left[L^p(T,\rho)\right]^*\equiv L^\pp(T,\rho^{1-\pp})$ under the $L^2(T,1)$-inner product.
A simple calculation shows that
$$
I^*_\mu(g)(\a)=\int_{\oS(\a)}gd\mu,
$$
so that (TreeCar) is equivalent to (\ref{milano}), and it suffices to consider positive $g$'s. By taking $g=\chi_{\oS(\a)}$, we see that (\ref{milano}) implies (TC). Set $\sigma(\a)=\rho^{1-\pp}(\a)\mu(\oS(\a))^\pp$. We rewrite (\ref{milano}) as
\begin{equation}
\label{melegnano}
\sum_{\a\in T}\left(\frac{1}{\mu(\oS(\a))}\int_{\oS(\a)}gd\mu\right)^\pp\sigma(\a)\le C(\mu)\int_{\oT}g^\pp d\mu
\end{equation}
Observe that the tree condition (TC) assumes the simple form
\begin{equation}
\label{dresano}
\sigma(\SS(\a))=\sum_{\b\ge\a}\sigma(\b)\le C(\mu)\mu(\oS(\a)).
\end{equation}
\begin{theorem}
\label{arbremagic}
Inequality (\ref{melegnano}) holds if and ony if condition (\ref{dresano}) holds. 

More generally, consider the maximal function $M$ defined on positive functions $g:\oT\to\RR$,
\begin{equation}
\label{pantigliate}
Mg(\a)=\max_{o\le\b\le\a}\frac{1}{\mu(\oS(\b))}\int_{\oS(\b)}gd\mu,\ \a\in T.
\end{equation}
Then, the inequality
\begin{equation}
\label{peschiera}
\sum_{\a\in T}\left(Mg(\a)\right)^\pp\sigma(\a)\le C(\mu)\int_{\oT}g^\pp d\mu
\end{equation}
is equivalent to (\ref{dresano}).
\end{theorem}
See \cite{ARS3} for some variations on this theme.

\begin{proof}
Condition (\ref{dresano}) is already necessary for the weaker (\ref{melegnano}). Let us turn to sufficiency.
The sublinear operator $M$ is bounded with unitary norm from $L^\infty(T,\mu)$ to $L^\infty(T,\sigma)$. By Marcinkievic interpolation, it suffices to show that $M$ is of weak type $1-1$. Let $\lambda>0$ and consider 
$E(\lambda)=\{\a\in T:\ Mg(\a)>\lambda\}$ and let $\Gamma$ be set of the minimal points (with respect to the partial order of $T$). Due to the tree structure,
$$
E(\lambda)=\bigsqcup_{\gamma\in\Gamma}\SS(\gamma).
$$ 
Hence,
\begin{eqnarray*}
\sigma(E(\lambda))&=&\sum_{\gamma\in\Gamma}\sigma(\SS(\gamma))\crcr
&\le&C(\mu)\sum_{\gamma\in\Gamma}\mu(\oS(\gamma))\crcr
&\le&C(\mu)\lambda^{-1}\sum_{\gamma\in\Gamma}\int_{\oS(\gamma)}gd\mu\crcr
&\le&C(\mu)\lambda^{-1}\int_{\oT(\gamma)}gd\mu
\end{eqnarray*}
\end{proof}
\section{Proof that (Car)$\implies$(TreeCar)}.
 Let $F,G$ be
holomorphic functions in $\mathbb{D}$,
\[
F(z)=\sum_{0}^{\infty}a_{n}z^{n},\ G(z)=\sum_{0}^{\infty}b_{n}z_{n}%
\]
Their Dirichlet inner product is
\[
\left\langle F,G\right\rangle _{\mathcal{D}}=\sum_{0}^{\infty}%
na_{n}\overline{b_{n}}=F(0)\overline{G(0)}+\int_{\mathbb{D}}F^{\prime}(z)\overline{G^{\prime}%
(z)}dA(z).
\]
The reproducing kernel of
$\mathcal{D}$ with respect to the product $\langle\cdot,\cdot\rangle
_{\mathcal{D}}$ is
\[
\phi_{z}(w)=1+\log{\frac{1}{1-w\bar{z}}}%
\]
i.e., if $f\in\mathcal{D}$, then
\[
f(z)=\langle f,\phi_{z}\rangle_{\mathcal{D}}=\int_{\mathbb{D}}f^\prime(w)\overline{\left(
1+\log{\frac{1}{1-\bar{z}w}}\right)^\prime}  dA(w)+f(0){.}%
\]
\begin{lemma}
\label{lemmaRap}Let ${\rho}$ be an admissible weight, $1<p<\infty$. If $G\in B_{p{^{\prime}}}({\rho}^{1-p{^{\prime}}})$, then
\begin{equation}
G(z)=\langle G,\phi_{z}\rangle_{\mathcal{D}}.\label{eqRap}%
\end{equation}
\end{lemma}

Now, let $\mu$ be a positive bounded measure on $\mathbb{D}$ and define
\[
\langle F,G\rangle_{\mu}=\langle F,G\rangle_{L^{2}(\mu)}=\int_{\oDD%
}F(z)\overline{G(z)}d\mu(z).
\]
The measure $\mu$ is Carleson for $B_{p}({\rho})$ if and only if
\[
Id:B_{p}(\rho)\rightarrow L^{p}(\mu)
\]
is bounded. In turn, this is equivalent to the boundedness, with the same
norm, of its adjoint $\Theta=Id^{\ast}$,
\[
\Theta:L^{{p^{\prime}}}(\mu)\rightarrow\left(  B_{p}({\rho})\right)  ^{\ast
}\equiv B_{p{^{\prime}}}({\rho}^{1-p{^{\prime}}})
\]
where we have used the duality pairings $\langle\cdot,\cdot\rangle
_{\mathcal{D}}$ and $\langle\cdot,\cdot\rangle_{\mu}$, and the assumption (dual) on $\rho$.

By {Lemma~\ref{lemmaRap}},
\begin{align*}
\Theta G(z)  & =\langle\Theta G,\phi_{z}\rangle_{\mathcal{D}}=\langle
G,\phi_{z}\rangle_{L^{2}(\mu)}\\
& =\int_{\oDD}\left(  1+\log{\frac{1}{1-z\bar{w}}}\right)  G(w)d\mu(w).
\end{align*}
Without loss of generality, assume that 
$\mbox{supp}(\mu)\subseteq\{z\colon|z|\ge1/2\}$. 
Consider, now, functions $g\in L^{p^\prime}(\mu)$, having the form
$$
g(w)=\frac{|w|}{\overline{w}}h(w)
$$ 
where $h\ge0$ and $h$ is constant on each box $\alpha\in T$, 
$h|_\alpha=h(\alpha)$. The boundedness of $\Theta$ implies
\begin{align*}
C\left(\int_{\oT}|h|^\pp d\mu\right)^{\frac{1}{p^\prime}}
&=\left(\int_{\oDD}|g|^{p^\prime}d\mu\right)^{\frac{1}{p^\prime}}\\
&\ge\|\Theta g\|_{B_{p^\prime}(\rho^{1-p^\prime})}\\
& \ge\left(\int_{\DD}\left|\int_{\oDD}\frac{1-|z|^2}{1-z\overline{w}}|w|h(w)d\mu(w)
\right|^{p^\prime}\rho(z)^{1-p^\prime}\frac{dA(z)}{(1-|z|^2)^2}
\right)^{\frac{1}{p^\prime}}\\
\end{align*}
For $z\in\mathbb{D}$, let $\alpha(z)\in T$ be the Whitney box
containing $z$. By elementary estimates,
\begin{equation}
\label{equseful}
\mbox{Re}\left(\frac{|w|(1-|z|^2)}{1-z\overline{w}}\right)\ge 0
\end{equation}
if $w\in\oDD$, and
$$
\mbox{Re}\left(\frac{|w|(1-|z|^2)}{1-z\overline{w}}\right)\ge 
c>0,\ \mbox{if\ }w\in \DoS(\alpha(z))
$$
for some universal constant $c$. 
Using this, and the fact that all our Whitney boxes have comparable
hyperbolic measure, $\int_\a\frac{dA(z)}{(1-|z|^2)^2}\sim1$, we can continue the chain of inequalities
\begin{align*}
& \ge \left(\int_{\mathbb
    D}\left|\int_{\DoS(\alpha(z))}\frac{1-|z|^2}{1-z\overline{w}}|w|h(w)d\mu(w)
\right|^{p^\prime}\rho^{1-p^\prime}(z)
\frac{dA(z)}{(1-|z|^2)^2}\right)^{\frac{1}{p^\prime}}\\
& \ge
c\left(\int_{\mathbb D}\left(\int_{
\oS(\alpha(z))}
h(\omega)d\mu(\omega)\right)^{p^\prime}\rho^{1-p^\prime}(z)
\frac{m(dz)}{(1-|z|^2)^2}
\right)^{\frac{1}{p^\prime}}\\
& \ge
c\left(\sum_{\alpha}\left(\int_{
\oS(\alpha(z))}
h(\omega)d\mu(\omega)\right)^{p^\prime}
\rho(\alpha)^{1-p^\prime}
\right)^{\frac{1}{p^\prime}}
\end{align*}
The chain of inequalities above shows that
$$
{I_\mu}^*\colon L^{p^\prime}(\mu)\rightarrow{L^{p^\prime}(\rho^{1-p^\prime})}
$$
is a bounded operator.
In turn,
this is equivalent to the boundedness of
$$
{\cal I}\colon L^p(\rho)\rightarrow L^p(\mu),
$$
which is (TreeCar).

This end the proof of Theorem \ref{main}.

\medskip

The duality argument used here ceases to work for the Hardy space.  If one could prove Theorem \ref{main} for the Hardy space, it would follow that the radial variation of a function $f$ in $H^2$ be finite $a.e.$ on $\partial\DD$, since the circular measure on $\partial\DD$ is Carleson for $H^2$. But there are functions in $H^2$ with infinite radial variation $a.e.$ on $\partial\DD$, so we have reached a contradiction.

See \cite{ARS1} and \cite{AR} for more information on this aspect and for explicir examples of Carleson measures for $H^2$ which do not satisfy (TC).

\end{document}